\documentclass[12pt,reqno]{amsart}
\usepackage{amsmath} \usepackage{amsfonts} \makeindex
\newcommand{\bd}{\partial}

\newcommand{\mc}{\mathcal}

\newcommand{\ol}{\overline}

\newcommand{\sphere}{\mathbb{C}_\infty}

\def\int{\mathbb{Z}}

\newcommand{\e}{\varepsilon}
\newcommand{\ucirc}{\mathbb{S}^1}
\newcommand{\jump}{\mathrm{JUMP}}
\newcommand{\jm}{\mathrm{JM}}
\newcommand{\0}{\emptyset}
\newcommand{\bbc}{{\mathbb C}}
\newcommand{\bbd}{{\mathbb D}}

\newcommand{\ch}{\mathrm{CH}}
\newcommand{\card}{\mathrm{card}}

\newcommand{\ord}{\mathrm{ord}}
\newcommand{\len}{\mathrm{len}}
\newtheorem{thm}{Theorem}[section]
\newtheorem{lem}[thm]{Lemma}
\newtheorem{cor}[thm]{Corollary} 
\newtheorem{prop}[thm]{Proposition}

\newtheorem{ass}[thm]{Standing Assumption}
\theoremstyle{definition}

\newtheorem{rem}[thm]{Remark}
\newtheorem{dfn}[thm]{Definition}

\begin{document}

\title
{Wandering Polygons and Recurrent Critical Leaves}
\author[D.K.~Childers]{Douglas K.~Childers}
\address[Douglas K.~Childers]
{Department of Mathematics\\ University of Alabama at Birmingham\\
Birmingham, AL 35294-1170} \email{childers@math.uab.edu}

\date{\today}

\begin{abstract}
Let $T$ be a finite subset of the complex unit circle $\ucirc$, and
define $f: \ucirc \mapsto \ucirc$ by $f(z) = z^d.$  Let $\ch (T)$
denote the convex hull of $T.$  If $\card(T) = N \geq 3,$ then
$\ch(T)$ defines a polygon with $N$ sides.  The $N$-gon $\ch(T)$ is
called a \emph{wandering $N$-gon} if for every two non-negative
integers $i \neq j,$ $\ch(f^i(T))$ and $\ch(f^j(T))$ are disjoint
$N$-gons.

A non-degenerate chord of $\ucirc$ is said to be \emph{critical} if
its two endpoints have the same image under $f.$  Then for a
critical chord, it is natural to define its (forward) orbit by the
forward iterates of the endpoints. Similarly, call a critical chord
\emph{recurrent} if one of its endpoints is recurrent under $f.$ The
main result of our study is that a wandering $N$-gon has at least
$N-1$ recurrent critical chords in its limit set (defined in a
natural way) having pairwise disjoint, infinite orbits.

Using this result, we are able to strengthen results of Blokh, Kiwi
and Levin about wandering polygons of laminations. We also discuss
some applications to the dynamics of polynomials.  In particular,
our study implies that if $v$ is a wandering non-precritical vertex
of a locally connected polynomial Julia set, then there exists at
least $\ord(v) -1$ recurrent critical points with pairwise disjoint
orbits, all having the same $\omega$-limit set as $v.$  Thus, we
likewise strengthen results about wandering vertices of polynomial
Julia sets.
\end{abstract}

\maketitle

\section{Introduction} \label{intro}
It is well-known that the dynamical behavior of critical points is
closely related to a number of interesting dynamical phenomena
exhibited by rational functions. One of the first observed examples
is that the immediate basin of an attracting cycle always contains a
critical point \cite{Fatou}. In the same paper, Fatou also proved a
result that implies that every indifferent cycle has an associated
critical point. Moreover, by the work of Douady \cite{Bl} and
Shishikura \cite{Sh}, the degree of a polynomial bounds the total
number of Cremer cycles and cycles of bounded Fatou components. A
nice proof of this inequality for polynomials with connected Julia
sets was given by Kiwi in \cite{K2}.  He obtained the inequality by
showing that the number of critical points bounds the total number
of attracting and indifferent cycles.  Within the study of rational
dynamical systems, another area where the importance of critical
points can be seen is in the investigation of the typical behavior
of points (see, for example, \cite{L1}, \cite{GPS}, \cite{BMO1},
\cite{BMO2}, \cite{BM}).

More specifically, the recurrence of critical points has recently
been recognized as an important dynamical property. For example,
Ma\~{n}\'{e} proved that a rational map is backward stable at every
non-parabolic point which is not in the $\omega$-limit set of a
recurrent critical point \cite{Ma}. This result proves relevant to
the topology and dynamics of Julia sets. In particular, Ma\~{n}\'{e}
showed that his result on backward stability implies that the
boundary of a Siegel disk (or respectively a Cremer point) is always
contained in the $\omega$-limit set of some recurrent critical
point. Thus, the presence of \emph{recurrent} critical points is
necessary for some types of complicated dynamics to occur for
rational maps. Other results in which the behavior of recurrent
critical points is shown to play an important role in dynamics can
be found in \cite{CJY}, \cite{BMO2} or \cite{BM}.

For the most part, these examples deal with the dynamics of rational
functions from an analytic point of view. However, critical points
and their analogs have also been shown to play a significant role in
the combinatorial and topological aspects of rational dynamical
systems. The aim of this paper is to establish and specify their
impact on the existence and behavior of wandering polygons in the
complex unit circle, an important phenomenon in the combinatorial
theory of polynomials.

In what follows, we use standard terminology without a formal
introduction, while defining the less standard notions. Always let
$f:\ucirc \mapsto \ucirc$ denote a map on the complex unit circle,
$\ucirc$, defined by $f(z) = z^d$ ($d \geq 2$). Thus $d$ will always
be the degree of $f.$  Also, for a set $B \subset \ucirc,$ let
$\ch(B)$ denote the convex hull of $B,$ and define $B_i = f^i(B).$
Finally, two closed sets $A,B \subset \ucirc$ are said to be
\emph{unlinked} if $\ch(A) \cap \ch(B) = \0$.

\begin{dfn}[Wandering Polygon] \label{def wand N-gon}
Let $T \subset \ucirc$ be a finite set with $N = \card(T) \geq 3.$
Then $T$ \emph{forms a wandering $N$-gon} if

\begin{enumerate}

\item \textbf{Non-Precritical}: For every $i$, $\card(T_i) = N.$

\item \textbf{Pairwise Unlinked}: For every $i \neq j$, $T_i$ and $T_j$ are
unlinked.

\end{enumerate}
\end{dfn}

Thurston proved that wandering $N$-gons do not exist under the map
$z \mapsto z^2$ \cite{Th}. About this result, he writes, ``This
theorem is invaluable for classifying laminations in the degree-two
case. The question whether this generalizes to higher degree, and if
not what is the nature of counterexamples, seems to me the key
unresolved question about invariant laminations'' \cite{Th}. After
Thurston, Kiwi was the first to make significant progress towards
understanding the necessary conditions for the presence of a
wandering polygon. This was done in the language of polynomial Julia
sets \cite{K1}, although stated below is his result in terms of a
wandering polygon (compare Theorem~A.2 in \cite{K3}).

\begin{thm}[Kiwi~\cite{K3}]\label{wandboundKiwi}
If $T \subset \ucirc$ forms a wandering polygon under $f$, then
$\card(T) \leq d.$
\end{thm}

Other results in this direction were obtained by Levin in \cite{L1}
(who was the first to discover a more general setting where
wandering polygons do not exist, even when $d \geq 3$), Blokh and
Levin in \cite{BL1} and \cite{BL2} (where several co-existing
polygons are considered), and Blokh in \cite{B1} (where the cubic
case is considered in detail). We discuss these results while
observing applications of our results in Section~\ref{Sec Applic},
where full definitions of necessary notions are given (including
that of a \emph{lamination}).

More recently, Blokh and Oversteegen have shown the existence of
wandering polygons under the maps $f(z) = z^d$ of degree $d \geq 3$
\cite{BO1}; as well as, the existence of an uncountable family of
polynomials with \emph{dendritic} (i.e. tree-like locally connected)
Julia sets which admit wandering non-precritical branch points
(equivalently, one can say that the corresponding laminations admit
wandering polygons). This justifies our interest in the dynamics of
wandering polygons. To study their dynamics, we use the notion of a
polygon's limit set.  From now on, let $T \subset \ucirc$ form a
wandering $N$-gon.

\begin{dfn}[Limit Set]\label{Def LimitSet}
Define the \emph{limit set} of $\ch(T)$, denoted by
$\mc{L}_{\omega(T)}$, as the collection of chords (including
degenerate chords - i.e. points in $\ucirc$) of $\ucirc$ such that
each element of $\mc{L}_{\omega(T)}$ is the limit set of chords from
a subsequence of $\ch(T_i).$ Any element of $\mc{L}_{\omega(T)}$ is
called a \emph{leaf}.
\end{dfn}
A non-degenerate chord of $\ucirc$ is called \emph{critical} if both
of its endpoints have the same image under $f.$  For a critical
chord $q,$ define its \emph{critical value} as the unique image of
the endpoints of $q,$ denoted $f(q).$  Critical chords in the limit
set $\mc{L}_{\omega(T)}$, called \emph{critical leaves}, play a
significant role in the dynamics of $T.$ Some of the intimate
relationship between critical chords and wandering polygons is
revealed in \cite{BL1}, \cite{K3}, and \cite{Th}. We further specify
this relationship by proving the existence of $N-1$ \emph{recurrent}
critical leaves, which are \emph{recurrent} critical chords in the
limit set $\mc{L}_{\omega(T)}.$  Below, we give the definition of a
recurrent critical chord.

\begin{dfn}[Dynamics of a Critical Chord]
Let $\ell$ be a critical chord with critical value $f(\ell) = v.$
Define the $\omega$-limit set of $\ell$ by $\omega(\ell) =
\omega(v).$ We say that $\ell$ is \emph{recurrent} if $\omega(\ell)
\cap \ell \neq \0.$
\end{dfn}

As a result of our study, we can state Theorem~\ref{MTheorem1},
below, which follows by putting together Theorem~\ref{(N-1)critical}
and the Recurrence Theorem~\ref{Recurrence Theorem}.

\begin{thm}\label{MTheorem1}
If $T \subset \ucirc$ forms a wandering $N$-gon, then there exist at
least $N-1$ recurrent critical leaves with pairwise disjoint
infinite orbits, all having the same $\omega$-limit set $X \subset
\ucirc.$  Moreover, $p \cap X \neq \0$ for every $p \in
\mc{L}_{\omega(T)}.$
\end{thm}

More generally, we also consider the possibility of several
wandering polygons co-existing in the following way.

\begin{dfn}
Let $\Gamma$ be a collection of pairwise disjoint finite subsets of
$\ucirc.$  Then $\Gamma$ is said to \emph{form a wandering
collection of polygons} if
\begin{enumerate}

\item
For every set $T \in \Gamma,$ $T$ forms a wandering polygon under
$f.$

\item
For every two different sets $A, B \in \Gamma$ and any two integers
$n,m \geq 0$, $A_n$ and $B_m$ are unlinked.
\end{enumerate}
\end{dfn}

Given a wandering collection $\Gamma$ of polygons, define
$\mc{L}_{\omega(\Gamma)} = \cup_{T \in \Gamma} \mc{L}_{\omega(T)}.$
Let $R$ denote the maximum number of recurrent critical chords in
$\mc{L}_{\omega(\Gamma)}$ with pairwise disjoint infinite orbits.
Let $\Omega$ be the collection of distinct $\omega$-limit sets of
recurrent critical chords in $\mc{L}_{\omega(\Gamma)}.$  Using
Theorem~\ref{MTheorem1} and a counting argument, we prove in
Section~\ref{SevPolysSec} the following.

\begin{thm}\label{MTheorem2}
If $\Gamma$ is a wandering collection of polygons, then
$$\card(\Gamma) \le \sum_{T \in \Gamma} (\card(T) - 2) \leq R -
\card(\Omega) \leq d - 1 -\card(\Omega).$$
\end{thm}
Our results allow us to improve some known results on wandering
polygons in laminations \cite{B1}, \cite{BL1}, \cite{BL2},
\cite{K1}, \cite{K3}, \cite{L1}, \cite{Th}. Consequently, we also
improve results on wandering vertices in locally connected
polynomial Julia sets. The history of this subject and applications
of our results are discussed in Section~\ref{Sec Applic}, in the
framework of laminations.

\subsection{Acknowledgements}
I wish to thank the UAB Laminations Seminar, including Dr.~Alexander
Blokh, advisor Dr.~John Mayer, and Dr.~Lex Oversteegen, for their
helpful comments on an earlier version of this paper.  In addition,
I would particularly like to thank Dr.~Blokh and Dr.~Mayer for their
part in the editing.  Finally, I would like to express my sincere
appreciation to Dr.~Blokh.  His ideas were the inspiration for this
study, while his suggestions provided the motivation, and his
observations improved the conclusion.
\subsection{Applications: Invariant Laminations and Polynomial
Julia Sets}\label{Sec Applic}

As previously mentioned, the motivation for studying wandering
polygons lies on their implications in polynomial dynamics. To
explain this connection, one needs to introduce the notion of an
\emph{invariant lamination}. In our study, we will not define an
invariant lamination as Thurston did in \cite{Th}. Instead, we will
define laminations as certain equivalence relations of $\ucirc$
(compare \cite{Do}, \cite{McM}). Recall that two closed sets $A, B
\subset \ucirc$ are said to be \emph{unlinked} if $\ch(A) \cap
\ch(B) = \0.$

\begin{dfn}[Invariant Lamination]\label{Def Lam}
An equivalence relation $\sim$ of $\ucirc$ is called an invariant
lamination if it satisfies the following four conditions.

\begin{enumerate}
\item

\textbf{Closed:} The graph of $\sim$ is closed in $\ucirc \times
\ucirc.$

\item

\textbf{Pairwise Unlinked:} Any two equivalence classes are
unlinked.

\item

\textbf{Forward Invariant:} If $B \subset \ucirc$ is an equivalence
class then $f(B)$ is also an equivalence class.

\item

\textbf{Gap Invariant:}  If $B \subset \ucirc$ is an equivalence
class with $\card(B) \geq 3$, then $f|_B$ is a covering map with
positive orientation.

\end{enumerate}
\end{dfn}

Given an invariant lamination $\sim$, define the quotient map $\pi:
\ucirc \mapsto J_\sim=J$, where $J = \ucirc/\sim$ denotes the
quotient space defined by $\sim.$ By the closed property of $\sim$,
$J$ is a locally connected continuum (by a \emph{continuum} we mean
a compact, connected metric space). Without loss of generality,
since $\sim$ is unlinked, we can assume $J$ is imbedded in the
Riemann Sphere $\sphere.$ The forward invariance of $\sim$ allows us
to define the \emph{induced} map $F: J \mapsto J$ by $F = \pi \circ
f \circ \pi^{-1}.$ The dynamics of $F$ is quite similar to the
dynamics of polynomials on their Julia sets.  Thus, such a continuum
$J$ is a \emph{topological Julia set} (of $F$). Fittingly, each
component of $\sphere \setminus J$ is called a \emph{Fatou domain}.

It can be shown that every locally connected polynomial Julia set is
a topological Julia set. In fact, for the locally connected Julia
set $J_P$ of a polynomial $P$ of degree $d$, we can naturally define
its \emph{corresponding invariant lamination} in the following way.
Let $\bbd\subset \sphere$ be the open unit disk. Then infinity is a
super-attracting fixed point of $P$. Let us denote its basin of
attraction by $A_{\infty}$, and define the Riemann map $\psi: \bbd
\mapsto A_{\infty}$ such that $\psi^{-1} \circ P \circ \psi = z^d.$
Note that $\bd A_{\infty} = J_P.$ Since $J_P$ is locally connected,
$\psi$ extends to a continuous surjection $\hat{\psi}: \ol{\bbd}
\mapsto \ol{A}_{\infty}.$  The \emph{corresponding (to $P$)
invariant lamination} $\sim$ is defined by the equivalence classes
$\{\hat{\psi}^{^{-1}}(y)\}_{y \in J}$ (it is easy to see that $\sim$
has properties (1) - (4)). We can think of $\pi =
\hat{\psi}|_{\ucirc}$ as the quotient map $\pi: \ucirc \mapsto J_P
\approx J_\sim=\ucirc/\sim$ such that, on the Julia set $J,$ $P|_J =
\pi \circ f \circ \pi^{-1}$ (even though $\pi$ is not 1-to-1, the
formula on the right is well-defined).

In the context of laminations that are not defined by polynomials,
let $J_\sim=J\subset \sphere$ denote the topological Julia set of
the induced map $F.$ We will always use $\sim$ to denote an
invariant lamination. Given $J$, $F$, and the corresponding $\sim$,
use $\pi$ to denote the quotient map defined by $\sim$ such that $F
= \pi \circ f \circ \pi^{-1}.$ Let us give some basic definitions
for points in a topological Julia set.

For a point $x \in J,$ define the \emph{order} of $x$ by $\ord(x) =
\card(\pi^{-1}(x)).$ Note that $\ord(x) = \card(~Q: Q
\mathrm{~is~a~component~of~} J \setminus \{x\}~).$ A point $b \in J$
is called an \emph{endpoint} (of $J$) if $\ord(b)=1$, a
\emph{cutpoint} (of $J$) if $\ord(b)\geq 2$, and a \emph{vertex} or
\emph{branch point} (of $J$) if $\ord(b) \geq 3.$

\begin{dfn}[Wandering Vertex]
A point $b \in J$ is called a \emph{wandering vertex} if $b$ is a
vertex of $J$ such that for every non-negative integer $i$, $F^i(b)$
is neither critical nor periodic.
\end{dfn}
\begin{rem}
It follows that if $b$ is a wandering vertex, then $\ord(b) =
\ord(F^i(b))$ for every $i \geq 0.$
\end{rem}

Clearly, the $\sim$-classes that form wandering polygons correspond
to wandering vertices in $J;$ both of which are, necessarily,
non-precritical. Kiwi studied the dynamics of such vertices in
\cite{K1}. The main result of his study, stated below, provides us
with yet another example of the influence critical points have on
the dynamics of polynomials. We need the following definition to
state his result.
\begin{dfn}[Narrow]
A point $x \in J$ is \emph{narrow} if it does not lie in the
boundary of any bounded Fatou domain and $\ord(x) = 1.$
\end{dfn}

\begin{rem}
A point $x\in J$ is narrow iff for every arc $I\subset J$ with $x\in
I$, $x$ is an endpoint of $I$.
\end{rem}

\begin{thm}[Kiwi~\cite{K1}]\label{KiwiNarrowTh}
Let $J$ be the locally connected Julia set of a degree $d$
polynomial, and suppose $b$ is a wandering vertex of $J.$  Then
there exist at least $\ord(b)-1$ narrow critical values in the
$\omega$-limit set of $b.$  Consequently, $\ord(b) \leq d.$
\end{thm}

We previously mentioned that the presence of recurrent critical
points is necessary for some specific types of complicated dynamics
to occur in rational maps. It turns out that wandering vertices are
included in this class of dynamics. Indeed, the $\omega$-limit set
of a wandering vertex in $J$ coincides with the limit set of some
recurrent critical point \cite{BL2}. In fact, the $\omega$-limit set
of a wandering vertex coincides, moreover, with the $\omega$-limit
set of at least two recurrent critical points. This was shown by
Blokh in \cite{B1} for cubic topological Julia sets. We show this
for arbitrary degree. In addition, we are able to strengthen these
results, while simultaneously improving the bounds given in
Theorem~\ref{KiwiNarrowTh}.
\begin{thm}\label{MT2 Intro}
Let $J$ be a topological Julia set, and suppose $b$ is a wandering
vertex of $J.$ Then there exists, at least, $\ord(b)-1$ recurrent
critical points from pairwise disjoint grand orbits in $\omega(b),$
each one having a narrow critical value and with their
$\omega$-limit set coinciding with $\omega(b).$
\end{thm}

\noindent Theorem~\ref{MTheorem1} is the fundamental step in proving
Theorem~\ref{MT2 Intro}; although, it also relies upon
Proposition~\ref{Cor critleaf=narrow}.

It was Levin who introduced, in \cite{L1}, the idea of studying the
dynamics of wandering polygons that come from a lamination by means
of studying the dynamics of the map on its topological Julia set
(the main aim of \cite{L1} was to introduce and study an important
notion of \emph{backward stability}). Ideas of \cite{L1} were
further developed in \cite{BL1} and \cite{BL2}. In particular, this
was done in \cite{BL1} to consider, for the first time, the
possibility of several wandering polygons (with pairwise disjoint
orbits) co-existing in an invariant lamination. To state the main
result from \cite{BL1} in this direction, we need the following
definition.

\begin{dfn}[All-critical]\label{def allcrit}
Let $\sim$ be an invariant lamination.  A non-degenerate equivalence
class $C$ of $\sim$ is said to be \emph{all-critical} in $\sim$,
denoted \emph{all-c}, if all points in $C$ have the same image under
$f.$
\end{dfn}

All-critical classes correspond under the quotient map to critical
points in the topological Julia set whose values (i.e. images) are
endpoints. Given an invariant lamination $\sim,$ let $k_{\sim}$ be
the maximum number of all-c classes from pairwise disjoint infinite
orbits.

\begin{thm}[Blokh and Levin~\cite{BL1}]\label{BL's Ineq}
Let $\Gamma$ be a non-empty collection of equivalence classes of
$\sim$ forming wandering polygons with pairwise disjoint orbits.
Then
$$\card (\Gamma)\le \sum_{T \in \Gamma}(\card(T)-2) \leq
k_{\sim}-1 \le d-1.$$
\end{thm}

\noindent Theorem~\ref{BL's Ineq} strengthens Kiwi's results. The
approach taken by Blokh and Levin is very different than Kiwi's; it
involves \emph{growing trees}, which were introduced in \cite{L1}
and \cite{BL1}, and inspired by the theory of Hubbard trees (see
\cite{DH1}, \cite{DH2}). This approach has proven fruitful in the
further study of wandering polygons (see \cite{B1},\cite{BL2}).

As a result of our study, we can replace $k_{\sim}$ with a finer
characteristic of the lamination, which involves the number of
recurrent all-c classes and the number of their limit sets.  Thus,
let us introduce these dynamical notions.

\begin{dfn}[Recurrent]\label{Def:All-crit Recrnt}
Let $C$ be an all-c class in the invariant lamination $\sim.$ Define
its $\omega$-limit set by $\omega(C) = \omega(s)$ for any $s \in C.$
Then, $C$ is called \emph{recurrent} if $\omega(C) \cap C \neq \0.$
\end{dfn}

We will call an \textbf{infinite} subset $W \subset \ucirc$ a
\emph{recurrent all-c limit set}, with respect to the invariant
lamination $\sim$, if there is a recurrent all-c class $C$ in $\sim$
with $\omega(C) = W$ (it follows that the orbit of $C$ must be
infinite). Denote by $R_{\sim}$ the maximal number of recurrent
all-c classes in $\sim$ with pairwise disjoint infinite orbits; then
$R_{\sim}\le k_{\sim}.$ Now, let $\Omega_{\sim}$ be the collection
of \emph{all} the recurrent all-c limit sets.
Theorem~\ref{MTheorem2} implies Theorem~\ref{MT2: Ineq}, which
strengthens Theorem~\ref{BL's Ineq}.

\begin{thm}\label{MT2: Ineq}
Let $\Gamma$ be a collection of equivalence classes forming
wandering polygons with pairwise disjoint orbits. Then
$$\card (\Gamma)\le \sum_{T \in \Gamma} (\card(T)
-2)\le R_{\sim} - \card(\Omega_{\sim}) \le d - 1 -
\card(\Omega_{\sim}).$$
\end{thm}

Likewise, we strengthen the bound from \cite{BL1} on the number of
wandering vertices in a topological Julia set $J$. To this end, we
need some notation. Namely, let $R_J$ be the maximum number of
recurrent critical points in $J$ with narrow values and pairwise
disjoint infinite orbits. In addition, let $\Omega_J$ be the
collection of distinct $\omega$-limit sets of recurrent critical
points with narrow values and infinite orbits. As a corollary of
Theorem~\ref{MT2: Ineq}, we have the following.

\begin{thm}\label{MCOR2}
Let $J$ be a topological Julia set. Let $\Gamma$ be a non-empty
collection of wandering vertices in $J$ with pairwise disjoint
orbits. Then
$$\card (\Gamma) \le \sum_{b \in \Gamma} (ord(b) - 2) \le R_J -
\card(\Omega_J) \le d - 1 - \card(\Omega_J).$$
\end{thm}

\section{Holes in a Finite Set}\label{holefinSect}

In this section, we introduce some notation and make some
introductory observations. We take the measure of the full angle in
$\ucirc$ to be equal to $1$ and define the length of an arc $\gamma
\subset \ucirc$ as its angle measure, denoted $\len(\gamma)$. For
any two points $v, w \in \ucirc$, use $[v, w]$ to denote the arc in
$\ucirc$ running counterclockwise from $v$ to $w.$

For a closed set $B \subset \ucirc,$ call each component of $\ucirc
\setminus B$ a \emph{hole} in $B.$ Thus if two sets $Q$ and $B$ are
unlinked, then $Q$ is contained one of the holes in $B.$ The
following proposition makes this observation more precise.
\begin{prop} \label{holcontain}
Let $Q,B \subset \ucirc$ be unlinked sets such that
$\card(B),\card(Q) \geq 3.$ Suppose there exists a $\tau > 0$ such
that both $Q$ and $B$ have at least two holes of length $\geq \tau.$
Then $Q$ is contained in a hole in $B$ whose length is $> \tau.$
\end{prop}

If $B\subset S^1$ then we say that a map $g|_B$ is \emph{orientation
preserving} if it is injective and preserves the cyclic order on any
three points of $B.$  To continue, we introduce the following
notation, in Sections \ref{sss labhol} and \ref{Sec Im-Hol}, for a
finite set $B \subset \ucirc$ with $\card(B) = M \ge 2$.

\subsubsection{Ordering Holes by Size}\label{sss labhol}

Denote the holes in $B$ by $\{H_k(B)\}_{k=1}^{M}$ such that
$\len(H_k(B)) \geq \len(H_{k-1}(B))$ for every $k > 1$, and define
$s_k(B) = \len(H_k(B)).$ Even though the above numbering of the
$H_k$'s is not always well-defined, the numbers $s_k$ are; still, we
make a choice for the $H_k$'s, and fix it for the given set $B.$
Note that $\cup_{k=1}^M \overline{H_k(B)} = \ucirc$, so
$\sum_{k=1}^M s_k(B) = 1.$ Also, note that each hole $H_k(B)$ is an
open arc in the circle whose closure intersects $B$ at exactly two
points. These two points define a chord which is one of the edges of
the polygon $\ch(B).$ Let us denote this chord by $e_k(B).$  So the
chord $e_k(B)$ denotes the edge of the polygon $\ch(B)$ that
corresponds to the $k^{th}$-smallest hole in $B.$

\subsubsection{Image-Holes and Remainders} \label{Sec Im-Hol}

For every $k \in \{1,\ldots, M\}$, define $\tilde{s}_k(B) = min
\{s_k(B) - \frac{j}{d}\mid \frac{j}{d} \le s_k(B)$ for $j = 0,1,
\ldots , d-1\}.$ We call $\tilde{s}_k(B)$ the \emph{remainder of
$s_k(B)$}. Remainders are useful for the analysis of the lengths of
holes in $f(B).$ Next, for a hole $H=(u, w)$ in $B,$ the arc $(f(u),
f(w))$ is called the \emph{image-hole} of $H.$ This terminology is
justified in Lemma~\ref{Lem dRem=Size IM-H}, below, where we show
that the remainders, $\{\tilde{s}_k(B)\}_{k=1}^{M},$ of the lengths
of the holes in $B$ determine the lengths of the holes in $f(B).$ To
do this, first observe the following theorem (Theorem 2.1 in
\cite{M1}).

\begin{thm}\label{jamie} If $B\subset S^1$ is a closed set and
$f|_B$ is injective then $f|_B$ preserves orientation if and only if
$S^1\setminus B$ contains $d-1$ pairwise disjoint open intervals of
length $\frac{1}{d}.$
\end{thm}

Now we are ready to prove Lemma~\ref{Lem dRem=Size IM-H}.

\begin{lem}\label{Lem dRem=Size IM-H}
Let $B \subset \ucirc$ be a finite set with $\card(B) = M \geq 2$,
such that $f|_B$ is injective. If $f|_B$ is orientation preserving,
then for every $k \in \{1, \ldots, M\};$ the \emph{image-hole} of
$H_k(B)=(u, w)$ is a well-defined hole $(f(u), f(w))$ in $f(B)$ of
length $d s_k(B) \pmod 1 = d \tilde{s}_k(B).$ Moreover, $f|_B$ is
orientation preserving iff $\sum_{k=1}^{M} \tilde{s}_k(B) =
\frac{1}{d}.$
\end{lem}

\begin{proof} The first claim in the lemma is left to the reader.
Let us show that this claim implies the rest of the lemma. If $f|_B$
is orientation preserving then, by the claim, $\sum_{k=1}^{M}
\tilde{s}_k(B) = \frac{1}{d}$ because there is a 1-to-1
correspondence between the holes in $B$ and their image-holes, which
are the holes in $f(B)$ (note that the lengths of the holes in
$f(B)$ sum up to $1$).

Now, suppose that $\sum_{k=1}^M \tilde{s}_k(B) = \frac{1}{d}.$
Observe that $d(s_k(B) - \tilde{s}_k(B))$ is the maximum number of
pairwise disjoint open arcs of length $\frac{1}{d}$ that are
contained in $\cup_{k=1}^{M}H_k(B) = \ucirc \setminus B.$ So the
number of pairwise disjoint open arcs with length $\frac{1}{d}$ in
the complement of $B$ is $d[\sum_{k=1}^{M} (s_k(B) -
\tilde{s}_k(B))]=d[1 - \sum_{k=1}^{M}\tilde{s}_k(B)]=d-1$. By
Theorem~\ref{jamie}, this implies that $f|_B$ is orientation
preserving.
\end{proof}

The image-hole of $H$ does not necessarily coincide with the image
of $H$ under $f.$ However, $f(H)$ is equal to the image-hole of $H$
when $\len(H) < \frac{1}{d}.$ As an awkward consequence of the first
claim in Lemma~\ref{Lem dRem=Size IM-H}, the image-hole of the
$k^{th}$-smallest hole in $B$ is not necessarily the
$k^{th}$-smallest hole in $f(B).$

\subsection{Holes in a Wandering N-gon}

The motivation for studying holes comes from two properties of
wandering polygons. The first of which can easily be seen as a
consequence of the unit disk having finite area and is left to the
reader.

\begin{prop} \label{Prop wsizeto0}
If $T \subset \ucirc$ form a wandering $N$-gon then $s_k(T_i) \to 0$
as $i\to\infty$, for $1\le k \leq N-2$.
\end{prop}

The proposition below allows us to use the concept of image-hole to
introduce a well-defined bijection from the holes in $T_i$ onto the
holes in $T_{i+1}.$ Its proof relies upon Theorem~\ref{jamie}.

\begin{prop} \label{Prop wcycorder}
Let $\mc{A}$ be a collection of pairwise unlinked closed sets in
$\ucirc$ such that $f$ is injective on each set. Then there are no
more than $d-1$ sets $A \in \mc{A}$ for which $f|_A$ is not
orientation preserving.
\end{prop}

\begin{proof}
For any collection $\mc{A}$ of pairwise unlinked closed subsets of
$\ucirc$, let $c(\mc{A})$ denote the cardinality of the collection
of distinct sets $A \in \mc{A}$ for which $f|_A$ is not orientation
preserving. We want to show that for such every collection $\mc{A}$
such that $f$ is injective on each set, $c(\mc{A}) \leq d-1.$ We
prove this bound by means of induction on $d \geq 1.$

The base case is $d=1.$ If $f(z) = z,$ then, clearly, for every
collection $\mc{A}$ of pairwise unlinked sets, $c(\mc{A}) = 0.$ So
next, we assume by induction that for all $k < d,$ if $f(z) = z^k,$
then every collection $\mc{A}$ of pairwise unlinked sets such that
$f$ is injective on each set satisfies $c(\mc{A}) \leq k-1.$ Set
$f(z) =z^d$ and let $\mc{A}$ be such a collection for $f$. We want
to prove that $c(\mc{A}) \leq d-1.$

Without loss of generality, we suppose there exists a set $A \in
\mc{A}$ such that $f|_A$ is not orientation preserving. Note that by
assumption, $f|_A$ is injective.  Let $\mc{H}$ be the collection of
distinct holes in $A.$ For every $H \in \mc{H}$, define the
collection of pairwise unlinked sets $\mc{A}_H = \{B \in \mc{A} \mid
B \subset H\}.$ Recall that, since all sets in $\mc{A}$ are pairwise
unlinked, for every set $B \in \mc{A}$ with $B \neq A,$ there is a
unique $H \in \mc{H}$ such that $B \subset H.$ Consequently
$c(\mc{A}) -1  = \sum_{H \in \mc{H}} c(\mc{A}_H).$  We will show
that $\sum_{H \in \mc{H}} c(\mc{A}_H) \leq d-2.$

To do this, note that for every $H \in \mc{H},$ there is a unique
integer $j_H$ such that $\frac{j_H}{d} \leq \len(H) < \frac{j_H
+1}{d}.$ Since $f|_A$ is injective, while not orientation
preserving, we have that $\sum_{H \in \mc{H}} j_H \leq d-2$
(Theorem~\ref{jamie}). Thus it suffices to show that $c(\mc{A}_H)
\leq j_H$ for every $H \in \mc{H}.$ Pick $H \in \mc{H}$ and note
that (by definition of $j_H$) there must be a closed arc of length $
\frac{j_H + 1}{d}$ that contains $\ol{H}.$ Given such an arc
$\gamma$, we can define a continuous orientation preserving map
$\pi: \gamma \mapsto \ucirc$ that identifies only the endpoints of
$\gamma$ (in particular, so that $\pi|_H$ is injective) and defined
in such a way that $g \circ \pi = f|_{\gamma}$, where $g(z) = z^{j_H
+ 1}.$ Since $\pi|_H$ is injective and orientation preserving, the
collection of sets $\mc{V} = \{ \pi(B) \mid B \in \mc{A}_H\}$ is
pairwise unlinked and in a one-to-one correspondence (induced by
$\pi$) to the collection $\mc{A}_H.$ Moreover for every set $B \in
\mc{A}_H$, $\pi|_B$ is an orientation preserving bijection. Thus we
are done, by induction, because $j_H$ is an upper bound for the
number of distinct sets $V \in \mc{V}$ with $g|_V$ not being
orientation preserving.
\end{proof}

To continue in our study, Propositions~\ref{Prop wsizeto0} and
\ref{Prop wcycorder} allow us to make the following assumption
without loss of generality.

\begin{ass}\label{STASS 1}
If a set $T \subset \ucirc$ forms a wandering $N$-gon then we may
assume that $f|_{T_i}$ preserves the orientation and for every $k
\in \{1, \ldots, N-2\}$ that $\frac{1}{dN}>s_k(T_i)$, so that the
image-hole of $H_k(T_i)$ is of length $ds_k$ ($i=1, 2, \dots$).
\end{ass}

Let us use $T$ to denote a finite set which forms a wandering
$N$-gon and satisfies Standing Assumption~\ref{STASS 1}. Then the
image-holes of any hole in $T$ are all well-defined and not
degenerate for all powers of $f.$ Moreover, the lengths of the
image-holes of a hole $H$ depend only on the length of $H.$ Using
Standing Assumption~\ref{STASS 1}, we can show that the lengths of
any two distinct holes in $T_i$ are different.

\begin{prop}\label{Prop distinct sizes}
Let $T$ be as in Standing Assumption~\ref{STASS 1}.  For every
integer $i \geq 0$, the lengths of the holes in $T_i$ are pairwise
distinct; moreover, their remainders are also pairwise distinct. In
particular, $$ s_N(T_i) >s_{N-1}(T_i) > \ldots > s_1(T_i).$$
\end{prop}

\begin{proof}
By Lemma~\ref{Lem dRem=Size IM-H}, a hole in $T$ has infinitely many
image-holes of length greater than $\frac{1}{d}$. By Standing
Assumption~\ref{STASS 1}(2), any hole in $T_i$ whose length is
greater than $\frac{1}{d}$ is one of the two longest holes in $T_i$.
Hence, if two holes in $T$ have the same length, then there exist
two non-negative integers $q \neq p$ with $s_{N-1}(T_p) = s_N(T_p)
\geq s_{N-1}(T_q)= s_N(T_q)>\frac{1}{d},$ which contradicts the fact
that $T_p$ and $T_q$ are unlinked (being unlinked implies that
$H_{N-1}(T_p)$ or $H_N(T_p)$ is contained in a hole in $T_q,$ while
each hole in $T_q$ has length at most equal to $s_N(T_p) =
s_{N-1}(T_p)$). The claim about the remainders follows from
Lemma~\ref{Lem dRem=Size IM-H}.
\end{proof}

Using Proposition~\ref{Prop distinct sizes}, we will see that there
is a unique edge of $T_i$ that is `closest' to a critical chord
(among those not crossing the polygon).   To do this, we need the
following definitions.

\subsubsection{Critical Hole}

By Proposition~\ref{Prop distinct sizes}, there is a unique integer
denoted by $cr(i) = cr  \in \{1, \ldots , N\}$ such that
$\tilde{s}_{cr}(T_i)$ is minimal among all the remainders of holes
whose lengths are greater than $\frac{1}{d}.$  By Standing
Assumption~\ref{STASS 1}, we see that $cr \in \{N-1, N\}.$ We call
$H_{cr}(T_i)$ the \emph{critical hole} in $T_i.$ For certain
$T_i$'s, the image-hole of the critical hole will contain the
critical value of the critical leaf `closest' to the polygon
$\ch(T_i)$; later on, this will give us a lot of information on the
dynamics of this critical leaf.

Next, we introduce a notion of distance between two chords that are
disjoint in the unit open disk $\bbd.$  This notion may be used to
define a metric on any collection of chords that are pairwise
disjoint in $\bbd$, which is equivalent to the Hausdorff metric.

\subsubsection{The $\rho$-metric}

For two chords $p,q$ that are disjoint in $\bbd$, there exist a
unique component of $\ol{\bbd} \setminus (p \cup q)$ whose closure
contains $p \cup q.$  Denote the closure of this component by
$C_{p,q}.$ Note that $C_{p,q} \cap \ucirc$ forms two arcs, or
possibly an arc union a point. Let $\rho(p,q)$ denote the sum of the
lengths of these two arcs. We will call $\rho(p,q)$ the \emph{amount
of arc between $p$ and $q$}. Let us observe some immediate
properties of $\rho.$

\begin{lem} \label{propofrho} Let $p$ be a chord and
$\{q_i\}_{i=1}^{\infty}$ a sequence of chords disjoint from $p.$

\begin{enumerate}

\item If $l$ is a chord such that $l$ separates $p$ from $q_1$ in
$\bbd$ then $\rho(p,l) + \rho(q_1,l) = \rho(p,q_1).$

\item If $q_1$ and $p$ are both critical chords then $\rho(p,q_1) \geq
\frac{1}{d}.$

\item The sequence $q_i$ converges onto $p$ iff $\rho(p,q_i)$
converges to zero.

\end{enumerate}
\end{lem}
By Lemma~\ref{propofrho}, we can use $\rho$ to determine if the
chords from a subsequence of $\ch(T_i)$ converge onto a given
critical chord.  To this end, we need Lemma~\ref{rhoeqrem}, which
describes the relationship between the remainder of a specific hole
in $T_i$ and the amount of arc between its corresponding edge and
certain critical chords. The proof of this lemma is left to the
reader; recall that $e_k(T_i)$ denotes the edge of the polygon
corresponding to the hole $H_k(T_i).$

\begin{lem} \label{rhoeqrem}
Suppose for some $k \in \{1, \ldots, N\}$ there is a $j \in \{1,
\ldots, d-1\}$ such that $\frac{j+1}{d} > s_k(T_i) > \frac{j}{d}.$
Let $[c,d]$ be an arc of length $\frac{j}{d}$ contained in
$H_k(T_i)$ and define the critical chord $q = \ol{cd}.$ Then
$\rho(e_k(T_i),q)=\tilde{s}_k(T_i),$ and $f(q)$ lies in the
image-hole of $H_k(T_i).$
\end{lem}

\subsection{Critical Leaves Exist}\label{Sec CritLeavExist}

Recall that a critical chord is called a \emph{critical leaf} if it
appears as the limit of chords from a subsequence of the polygons
$\ch(T_i).$  In this section, we will prove that critical leaves
exist.  Also, we include a detailed account of how the critical
leaves appear.  We believe this is the foundation for understanding
the relationship between the wandering polygon and the dynamics of
its critical leaves.  We begin with the following definition.
\begin{dfn}\label{DEF T-Narr}
A point $x\in S^1$ is called \emph{$T$-narrow} if there exists a
subsequence $T_{i_n}$ such that $\lim_{n \to \infty} \ch(T_{i_n}) =
x$ and $x \in \cup_{j=1}^{N-1}H_j(T_{i_n})$ for every $n.$
\end{dfn}
Proposition~\ref{Cor critleaf=narrow} immediately follows from
Lemma~\ref{propofrho}, the properties of $\rho$, and the
definitions; the proof is left to the reader.

\begin{prop}\label{Cor critleaf=narrow}
Every critical leaf has $T$-narrow image.  Consequently, for every
critical leaf $p \in \mc{L}_{\omega(T)},$ $ \ch(T_i) \cap p = \0$
for every $i.$
\end{prop}

To conclude this section, we describe the relationship between
critical remainders and critical leaves.

\begin{thm} \label{Th CritLeafExist}
If $T$ forms a wandering $N$-gon then critical leaves exist. In
particular, there exists a positive $\e_T < \frac{1}{dN}$ such that
for every non-negative integer $i$,  if $\tilde{s}_{cr}(T_i)< \e_T$
then there is a unique critical leaf $q$ such that
$$\rho(e_{cr}(T_i), q) = \tilde{s}_{cr}(T_i).$$
Consequently, in this case, the critical value $f(q)$ lies in the
image-hole of the critical hole $H_{cr}(T_i).$
\end{thm}

\begin{proof} By Proposition~\ref{Prop wsizeto0}, for at least $N-2$
holes $H$ in $T$ the lower limit of the sequence of lengths of their
image-holes is $0.$ Fix one such hole $H.$ Denote by $H^i$ the
$i$-th image-hole of $H$ (so that $H^0=H$), and by $e^i$ the chord
connecting the endpoints of $H^i$. Then there exists a sequence
$k_i$ such that $\len(H^{k_i})=\min\{\len(H^j), 0\le j\le k_i\},$ so
$\len(H^{k_i})\searrow 0$. By Lemma~\ref{Lem dRem=Size IM-H}; the
remainders of the image-holes $H^{k_i-1}$ converge monotonically to
$0,$ while the lengths of the image-holes $H^{k_i-1}$ are greater
than $\frac{1}{d}.$ By refining this sequence and by using Lemma
\ref{propofrho} and the properties of $\rho$ (Lemma
\ref{propofrho}), we may assume that the chords $e^{l_i}$ converge
to a critical leaf $l$. Hence, critical leaves exist. Also observe
that since $T$ is wandering, any two critical leaves are disjoint in
$\bbd,$ so there are only finitely many of them.

By way of contradiction, suppose that the number $\e_T$ does not
exist. Since $\liminf \tilde{s}_{cr}(T_i) = 0,$ this means that
there exists a sequence $i_n$ such that $\tilde{s}_{cr}(T_{i_n})<
\frac{1}{n},$ while a unique critical leaf $q$ with
$\rho(e_{cr}(T_{i_n}), q) = \tilde{s}_{cr}(T_{i_n})$ does not exist.
By the properties of $\rho$ and because the $\rho$-distance between
distinct critical leaves is at least $\frac{1}{d};$ we see that if
$n$ is sufficiently large, and there is a critical leaf $q$ with
$\rho(e_{cr}(T_{i_n}), q) = \tilde{s}_{cr}(T_{i_n})<\frac{1}{n},$
then this leaf is unique. Thus for each $n,$ there is no critical
leaf $q$ with $\rho(e_{cr}(T_{i_n}), q) =
\tilde{s}_{cr}(T_{i_n})<\frac{1}{n}.$ However, as in the first
paragraph of the proof, we can now choose a subsequence along which
chords $e_{cr}(T_{i_n})$ will converge to a chord $q,$ which then
will have to be a critical leaf. Clearly, this contradiction implies
that the desired number $\e_T$ exists; we fix $\e_T$ and use in what
follows.
\end{proof}

\section{Jumping Critical Leaves}\label{JumpingHolSec}
We begin with extending Standing Assumption~\ref{STASS 1}, which
further requires that $s_{N-2}(T_i) < \e_T$ and can be done without
loss of generality, by Proposition~\ref{Prop wsizeto0}; from now on,
we assume that Standing Assumption~\ref{STASS ep-wand} holds for
$T.$
\begin{ass}\label{STASS ep-wand}
If a set $T \subset \ucirc$ forms a wandering $N$-gon then we may
assume that $f|_{T_i}$ preserves orientation and $\frac{1}{3dN}>
\e_T > s_{N-2}(T_i)$, so that the image-hole of $H_k(T_i), 1\le k\le
N-2$, is of length $ds_k$ ($i=1, 2, \dots$).
\end{ass}

The main objective of this section is to understand the dynamics of
\emph{jumping critical leaves}.  Jumping critical leaves are defined
later; at this point we introduce and study a few related notions.

\subsection{Set of Jumps} \label{jumpset}

The iterates $i$ for which $d s_{N-2}(T_i)> s_{N-2}(T_{i+1})$ are
called the \emph{jumps}; the set of jumps is denoted by
$\jump(T)=\{\jm_1(T), \jm_2(T), \dots\}.$ In fact, jumps are exactly
the times $i$ when Theorem~\ref{Th CritLeafExist} is satisfied (i.e.
$\tilde{s}_{cr}(T_i) < \e_T$) \emph{and} the image-hole of
$H_{cr}(T_i)$ is one of the $N-2$ smallest holes in $T_{i+1}$. To
study the set of jumps, we need to better describe the relationship
between the holes in $T_i$ and their image-holes in $T_{i+1}.$ The
following lemma gives a sufficient condition on the length
$s_{N-1}(T_i)$ for when $f(H_k(T_i)) = H_k(T_{i+1})$ for every $k
\in \{1, \ldots, N-1\}$ and, as a corollary, describes a useful
property of jumps.

\begin{lem} \label{small N-1}
Let $T$ be as in Standing Assumption~\ref{STASS ep-wand}.

\begin{enumerate}

\item If $\frac{1}{d^{m+1}N}>s_{N-1}(T_i)$ then for all $k\in\{1, \dots,
N-1\}$ and all $j\le m$ we have $f^j(H_k(T_i))=H_k(T_i+j)$.

\item $\lim_{i\to \infty} (\jm_{i+1}(T) - \jm_i(T))\to \infty$

\end{enumerate}

\end{lem}

\begin{proof} (1) First we show that if $\frac{1}{dN} > s_{N-1}(T_i)$
then for all $k \in \{1, \ldots, N - 1\}$ we have $f(H_k(T_i)) =
H_k(T_{i+1})$ and $s_k(T_{i+1}) = d s_k(T_i).$ Clearly, if
$s_{N-1}(T_i) < \frac{1}{dN}$ then $H_N(T_i)$ is the critical hole.
More importantly, $\sum_{k=1}^{N-1} \tilde{s}_k(T_i) <
\frac{N-1}{dN}$ and hence $\tilde{s}_N(T_i) = \frac{1}{d} -
\sum_{k=1}^{N-1} \tilde{s}_k(T_i) > \frac{1}{dN}.$ So $d s_N(T_i)
\pmod 1 = d \tilde{s}_N(T_i) > d s_{N-1}(T_i),$ which implies that
the $N-1$ smallest holes in $T_{i+1}$ have lengths $ds_{N-1}(T_i) >
ds_{N-2}(T_i) > \ldots > ds_1(T_i).$ Since by Proposition~\ref{Prop
distinct sizes}, the holes in $T_{i+1}$ have pairwise distinct
lengths; thus Lemma~\ref{Lem dRem=Size IM-H} completes the argument.
Claim (1) of the lemma now easily follows from the proven above.

Observe that by the definition of a jump, none of the moments $i,
i+1, \dots, i+m$ is a jump.

(2) Suppose that $n$ is very big. Then by Proposition~\ref{Prop
wsizeto0}, the holes $H_1(T_n), \dots, H_{N-2}(T_n)$ are very small.
Hence by Lemma~\ref{Lem dRem=Size IM-H}, the order among these holes
and their image-holes is the same, and the size, $ds_{N-2}(T_n)$, of
the image-hole of $H_{N-2}(T_n)$ is still very small. If $n$ is a
jump then $s_{N-2}(T_{n+1})<ds_{N-2}(T_n)$, which implies that in
addition to the $N-3$ image-holes of $H_1(T_n), \dots, H_{N-3}(T_n)$
(which are all shorter than $ds_{N-2}(T_n)$) there is \emph{at
least} one more hole in $T_{n+1}$ which is shorter than
$ds_{N-2}(T_n)$. However by Lemma~\ref{Lem dRem=Size IM-H}, the sum
of remainders of all the holes is $\frac{1}{d}$. Hence, there must
exist \emph{exactly} one more hole in $T_{n+1}$ which is shorter
than $ds_{N-2}(T_n)$; and so $s_{N-1}(T_{n+1})=ds_{N-2}(T_n),$ which
is very small; using this along with (1) (see the observation in the
end of the proof of (1)), there is a long stretch of numbers $n+1,
\dots$ that contains no jump, as desired.
\end{proof}

Using a similar argument, we describe in Proposition~\ref{Prop
Im-Hol(N-2)} the possibilities for the image-hole of each one of the
$N-2$ smallest holes in $T_i.$

\begin{prop}\label{Prop Im-Hol(N-2)}
For every $k \in \{1, \ldots, N-2\}$, exactly one of the following
must occur.
\begin{enumerate}

\item

$s_k(T_i)> \tilde{s}_{cr}(T_i)$ and $f(H_k(T_i)) =
H_{k+1}(T_{i+1}).$

\item

$\tilde{s}_{cr}(T_i) > s_k(T_i)$ and $f(H_k(T_i)) = H_k(T_{i+1}).$

\end{enumerate}

Furthermore, if $k \in \{1, \ldots, N-2\}$ is such that (1) holds,
then the image-hole of the critical hole is one of the $k$ smallest
holes in $T_{i+1}$ and $\tilde{s}_{cr}(T_i)<\e_T.$
\end{prop}

\begin{proof}
Let us introduce the following notation: (a) if
$s_{N-1}(T_i)<\frac{1}{d}$ then set $H'=H_{N-1}$ and $s'=s_{N-1}$;
(b) if $s_{N-1}(T_i)>\frac{1}{d}$ then set $H'$ to be the
\emph{non-critical} hole of length greater than $\frac{1}{d}$ and
$s'$ to be equal the remainder of this length. Let us show then, in
either case, that the length of the image-hole of $H'$ is
$ds'>s_{N-2}(T_{i+1}).$ Indeed, it is enough to consider the case
when the length of $H'$ is greater than $\frac{1}{d}.$ In this case

$$s'+\tilde{s}_{cr}(T_i)=1 -
\sum_{k=1}^{N-2} \tilde{s}_k(T_i)>1-\frac{N-2}{3dN}>\frac{1}{Nd},$$

\noindent which implies (by the choice of the critical hole) that
$s'>\frac{1}{2Nd},$ and so $s'>\tilde{s}_k(T_i)\,$ $(1\le k\le
N-2).$ In fact, in the case (b) (i.e., when $H'$ is the non-critical
hole of length grater than $\frac {1}{d}$), we always have that
$ds'=s_N(T_{i+1})$. However, in general (i.e., when either (a) or
(b) can take place), we have that only if (1) holds for $i$ is
$ds'=s_N(T_{i+1})$ (by (1) the image-hole of the critical hole of
$T_i$ is shorter than $ds'$).

Hence, the only way the image-holes of the $N-2$ smallest holes in
$T_i$ are not the $N-2$ smallest holes in $T_{i+1}$ is when (1)
holds; in which case, the image-hole of the critical hole is one of
the $k$ smallest holes in $T_{i+1},$ and so $\tilde{s}_{cr}(T_i) < d
\tilde{s}_{cr}(T_i) < \e_T.$
\end{proof}

By Proposition~\ref{Prop Im-Hol(N-2)}, the inequality $d
s_{N-2}(T_i) > s_{N-2}(T_{i+1})$ is necessary and sufficient for the
image-holes of the $N-2$ smallest holes in $T_i$ not to be the $N-2$
smallest holes in $T_{i+1}.$ Also, these times $i$ are
\emph{exactly} the times when the image-hole of the critical hole is
one of the $N-2$ smallest holes in $T_{i+1}.$ Note also that by
Proposition~\ref{Prop Im-Hol(N-2)} and Standing
Assumption~\ref{STASS ep-wand}, $\tilde{s}_{cr}<\e_T$ for any $i\in
\jump(T)$, so we can apply Theorem~\ref{Th CritLeafExist} to any
integer $i \in \jump(T).$ This implies the following theorem, whose
proof we leave to the reader.

\begin{thm} \label{TH jump-noj sumry}
Let $T$ be as in Standing Assumption~\ref{STASS ep-wand}.

\begin{enumerate}

\item
If $i \in \jump(T)$ then there exists a unique critical leaf $q_i$
such that for the chord $e_{cr}(T_i)$ we have
$\rho(e_{cr}(T_i),q_i)= \tilde{s}_{cr}(T_i)<\e_T.$  It follows that
the critical value $f(q_i)$ lies in the image-hole of the critical
hole, which is one of the $N-2$ smallest holes in $T_{i+1}.$

\item If $i \not \in \jump(T)$ then for every $k \in \{1,
\ldots, N-2\}$, $f(H_k(T_i)) = H_k(T_{i+1}).$

\end{enumerate}

\end{thm}

Below, we use the notation $q_i$ as in Theorem~\ref{TH jump-noj
sumry}. At the time $i$ of a jump, the polygon $\ch(T_i)$ is
`leaf-like' and very close to a critical leaf $q_i$, while
$\ch(T_{i+1})$ is `point-like' and $f(q_i)$ is `trapped' in one of
the $N-2$ smallest holes in $T_{i+1}.$  We will see the dynamical
importance of these moments later in this section. First, we state a
corollary of Theorem~\ref{TH jump-noj sumry} (the proof is left to
the reader).

\begin{cor}\label{crit-iterate}
Let $j \in \jump(T)$ and suppose the integer $m > j$ is chosen so
that $i \not \in \jump(T)$ for every integer $i$ with $j < i < m.$
Then $f^{m-j}(q_j)$ lies in one of the $N-2$ smallest holes in
$T_m.$
\end{cor}

\subsection{Jumping Critical Leaf} \label{N-1jumpcrit}

We are interested in the critical leaves $l$ such that $q_i = l$ for
infinitely many $i \in \jump(T)$, which we call the \textbf{jumping
critical leaves}. It should be noted that while every critical leaf
has its image coinciding with a narrow critical value, a critical
leaf may not be one of the jumping critical leaves. However, the
critical leaves with narrow critical values that Kiwi found in
\cite{K1} correspond to the images of jumping critical leaves.  This
follows from our proof of Theorem~\ref{(N-1)critical}, which is
inspired by Kiwi's arguments from \cite{K1}.
Theorem~\ref{(N-1)critical} says that there are at least $N-1$
jumping critical leaves with pairwise disjoint orbits. Note that
there might be more than $N-1$ jumping critical leaves (e.g.,
consider the case when $T$ is strictly contained in some finite set
$P$, where $P$ also forms a wandering polygon). We show in
Section~\ref{Sec jumrec} that every jumping critical leaf is
recurrent.

\begin{thm}\label{(N-1)critical}
If $T$ forms a wandering $N$-gon, then there exist at least $N-1$
jumping critical leaves with pairwise disjoint orbits.
\end{thm}

\begin{proof}
By Proposition~\ref{Prop wsizeto0}, for every positive integer $n$
there exists $M_n$ such that $\frac{1}{d^{2n}N} > s_{N-2}(T_i)$ for
all $i \geq M_n.$  Also by Proposition~\ref{Prop wsizeto0}, we can
find a strictly increasing sequence $\{i_n\}_{n=1}^{\infty}$ of
integers such that for every $n$, $d s_1(T_{i_n}) > s_1(T_{i_n
+1}).$  We assume without loss of generality that $i_n \geq M_n$ for
every $n.$

Since there are finitely many critical leaves, by choosing
subsequences of $i_n,$ we may also assume that the leaf $q_{i_n}$
(defined in Theorem~\ref{TH jump-noj sumry}) is the same critical
leaf $l_0$ for every $n.$  Note that by Proposition~\ref{Prop
Im-Hol(N-2)}, all of this implies that $f(l_0) \in H_1(T_{i_n + 1})$
and $s_{N-1}(T_{i_n +1}) = d s_{N-2}(T_{i_n}) < \frac{1}{d^{2n -1}
N}.$ Then by applying Lemma~\ref{small N-1} to $T_{i_n + 1}$, we get
that $f^n(l_0) \in H_1(T_{i_n + n}).$

Set $\tau_n = s_1(T_{i_n + n}).$ By our assumptions and
Lemma~\ref{small N-1} we have $\tau_n<\frac{1}{d^nN}.$ The number
$\tau_n$ serves as a threshold below which the lengths $s_1(r),
\dots, s_{N-2}(r)$ of holes in $T_r$ must stay from some time on by
Proposition~\ref{Prop wsizeto0}. Choose $j_n(k)$ as the moment when
the size of the hole $s_k(j_n(k))$ drops below $\tau_n$ for the last
time (so that for $r>j_n(k)$ we have $s_k(r)<\tau_n$). We now show
that the finite sequence $j_n(1), \dots, j_n(N-2)$ is strictly
increasing.

In other words, we show that for every $n$ there exists a strictly
increasing finite sequence of integers $\{j_n(k)\}_{k = 1}^{N-2}$,
with $j_n(1) \geq i_n + n$, such that for each $k$,

\begin{enumerate}

\item[(a)] $s_k(T_{j_n(k)}) \geq \tau_n$, and

\item[(b)] for every $i > j_n(k)$, $\tau_n > s_k(T_i).$

\end{enumerate}

To this end, we fix $n$ and use finite induction on $k \in \{1,
\ldots, N-2\}.$ The base case is when $k=1.$ Take $j_n(1)$ to be the
maximal integer $i \geq i_n + n$ such that $s_1(T_i) \geq \tau_n,$
and note that $j_n(1)$ exists by Proposition~\ref{Prop wsizeto0}. By
our choice of $j_n(1)$, we have that $\tau_n > s_1(T_i)$ for all $i
> j_n(1).$  Moreover, we also have that $s_1(T_{j_n(1)}) \geq \tau_n
> s_1(T_{j_n(1) + 1})$, thus concluding the base case.

Now we assume that for some $k \in \{1, \ldots, N-3\}$, we have
found the integers $j_n(k)$ satisfying (a) and (b) above. Then by
Proposition~\ref{Prop Im-Hol(N-2)}~(1), we have that
$s_{k+1}(T_{j_n(k) +1}) = d s_k(T_{j_n(k)}).$ Since $s_k(T_{j_n(k)})
\geq \tau_n$, this implies $s_{k+1}(T_{j_n(k) + 1})> \tau_n.$  The
rest of the argument is analogous to the proof in the base case.

We can do this for every $n.$  Note that by Theorem~\ref{TH jump-noj
sumry}~(2), $j_n(k) \in \jump(T)$ for every $n$ and $k.$  Since the
number of critical leaves is finite, there is a collection of
critical leaves $\{l_1, \ldots, l_{N-2}\}$ and an infinite set $V
\subset \mathbb{N}$ such that if $n \in V$ then for every $k$,
$q_{j_n(k)} = l_k.$  Recall that $l_0$ is the critical leaf
corresponding to the sequence $i_n.$  We want to show that the
jumping critical leaves (given above) $l_0, \ldots, l_{N-2}$ have
pairwise disjoint orbits.  To see this, assume by way of
contradiction that there exists two non-negative integers $b,c$ such
that $f^b(l_{k_0}) = f^c(l_{k_1})$, where $k_0 \in \{0, \ldots,
N-2\}$, $k_1 \in \{1, \ldots, N-2\}$ and $k_0 \neq k_1.$

Since $V$ is infinite, we can fix $n \in V$ with $n >
\mathrm{max}\{b, c\}.$ Then for each $k \in \{1, \ldots, N-2\}$,
$q_{j_n(k)} = l_k.$  Also, Theorem~\ref{TH jump-noj sumry}~(1)
implies that $f(l_k)$ lies in the image-hole of the critical hole
$H_{cr}(T_{j_n(k)}).$  By the choice of our sequences and
Proposition~\ref{Prop Im-Hol(N-2)}, we get the inequality
$s_{N-1}(T_{j_n(k) + 1}) = d s_{N-2}(T_{j_n(k)}) > \tau_n
> s_k(T_{j_n(k) + 1}) \geq d \tilde{s}_{cr}(T_{j_n(k)}).$ Moreover
since $j_n(k) > i_n \geq M_n$, we have that $\frac{1}{d^{2n-1}N} >
s_{N-1}(T_{j_n(k) + 1}).$  By these inequalities, Lemma~\ref{small
N-1} and property~(b) imply the following.

For every $k \in \{1, \ldots, N-2\}$ and every positive integer $r <
2n$:
\begin{enumerate}
\item There are at least two holes in $T_{j_n(k) + r}$ of length
$\geq \tau_n$, while

\item $f^r(l_k)$ lies in a hole in $T_{j_n(k) + r}$ that has
length no greater than $\tau_n$, and

\item if $k < N-2$ then $j_n(k+1) > j_n(k) + 2n-1.$
\end{enumerate}

Define $\Gamma = \{T_{i_n + n}, T_{j_n(1) + 1}, \ldots, T_{j_n(1) +
2n -1}, \ldots, T_{j_n(N-2) + 1}, \ldots, T_{j_n(N-2) + 2n -1}\}.$
Note that the iterates in $\Gamma$ are pairwise distinct by
property~(3) and because $j_n(1) \geq i_n +n.$  In particular since
$T$ is a wandering $N$-gon, $\Gamma$ is a collection of pairwise
unlinked distinct sets. Thus since $\tau_n = s_1(T_{i_n + n}),$ and
since $T_{j_n(k) + r}$ satisfies properties~(1),(3) for every $k \in
\{1, \ldots, N-2\}$ and $r< 2n$, any two different sets $Q,B \in
\Gamma$ satisfies the conditions of Proposition~\ref{holcontain}
with $\tau = \tau_n.$  Hence, for any such two sets, $Q$ is
contained in a hole in $B$ whose size is bigger than $\tau_n.$

Recall that $n > \mathrm{max}\{b, c\}$.  It follows that $n \geq n-b
> 0$, $0 < r := n-b+c < 2n,$ and $f^n(l_{k_0}) = f^{n-b} \circ f^b
(l_{k_0}) = f^{n-b} \circ f^c(l_{k_1}) = f^r(l_{k_1}).$  Since $k_0
\in \{0, \ldots, N-2\},$ $B := T_{j_n(k_0) + n} \in \Gamma$ (where
$j_n(0) = i_n$).  Also, since $0< r < 2n$ and $k_1 \in \{1, \ldots,
N-2\}$, we have that $Q:= T_{j_n(k_1)+r} \in \Gamma.$  Observe that
by property~(3), $Q \neq B$ because $k_0 \neq k_1;$ so that $Q$ is
contained in a hole in $B$ whose size is bigger than $\tau.$ By the
symmetry of this argument, we also have that $B$ is contained in a
hole in $Q$ whose size is bigger than $\tau_n$.  However this
contradicts property~(2) of $Q$ (or $B$), since $f^n(l_{k_0}) =
f^r(l_{k_1}).$ \end{proof}

\subsection{Jumping Critical Leaves Are Recurrent}\label{Sec jumrec}

Recall that the $\omega$-limit set of a critical leaf $l$ is defined
by $\omega(l) = \omega(v)$, where $v = f(l) \in \ucirc$ is the
critical value.  We call a critical leaf \emph{recurrent} if
$\omega(l) \cap l \neq \0.$  Lemma~\ref{trapped seq} shows the
dynamical significance of critical values being `trapped' in one of
the $N-2$ smallest holes.  The proof of Lemma~\ref{trapped seq} is
left to the reader.

\begin{lem} \label{trapped seq}
Let $T$ be a wandering $N$-gon. Suppose $\{i_n\}_{n=0}^{\infty}$ is
a increasing sequence of non-negative integers such that
$\ch(T_{i_n})$ converges onto a leaf $p \in \ucirc.$ Let
$\{x_n\}_{n=0}^{\infty}$ be a sequence of points from $\ucirc$ such
that for every $n$, $x_n$ is in one of the $N-2$ smallest holes in
$T_{i_n}.$ Then the limit points of $\{x_n\}$ form a non-empty
subset of $p \cap \ucirc.$
\end{lem}

Now we are ready to prove the Recurrence Theorem.

\begin{thm} [Recurrence Theorem] \label{Recurrence Theorem}
All jumping critical leaves of a wandering $N$-gon are recurrent
with the same $\omega$-limit set.
\end{thm}

\begin{proof}
By Theorem~\ref{(N-1)critical} there are at least two jumping
critical leaves. Let $\mc{C}_{J}$ denote the collection of all
jumping critical leaves. Let $q \in \mc{C}_{J}$ and define the
collection $\mc{C}_{\omega(q)} = \{ l \in \mc{C}_{J} : l = q~
\mathrm{or}~ l \cap \omega(q) \neq \0 \}.$ Clearly, $\mc{C}_{J}$ is
non-empty, and $\omega(l) \subset \omega(q)$ for every $l \in
\mc{C}_{\omega(q)}.$ We show that $\mc{C}_{\omega(q)} = \mc{C}_{J},$
which concludes the argument since there are at least two jumping
critical leaves and since $q \in \mc{C}_{J}$ was arbitrary.

By definition $\mc{C}_{\omega(q)} \subset \mc{C}_J.$  Assume by way
of contradiction that $\mc{C}_{J} \setminus \mc{C}_{\omega(q)} \neq
\0.$  Then there exist infinitely many integers $i \in \jump(T)$
with $q_i \not \in \mc{C}_{\omega(q)}.$  Since $\mc{C}_{\omega(q)}
\neq \0$, we also know that there are infinitely many integers $i
\in \jump(T)$ with $q_i \in \mc{C}_{\omega(q)}.$ Hence we can find a
strictly increasing sequence $\{j_n\}_{n=1}^{\infty}$ of jumps such
that for every positive integer $n$, $q_{j_{2n-1}} \in
\mc{C}_{\omega(q)}$ while $q_{j_{2n}} \in \mc{C}_J \setminus
\mc{C}_{\omega(q)}$ and $i \not \in \jump(T)$ for all $i$ with
$j_{2n-1}< i < j_{2n}.$ Since there are only finitely many critical
leaves we may assume that for every $n$, $q_{j_{2n-1}}$ is the same
critical leaf $l \in \mc{C}_{\omega(q)}$ and $q_{j_{2n}}$ is the
same critical leaf $r \in \mc{C}_J \setminus \mc{C}_{\omega(q)}.$

For each $n$, Theorem~\ref{TH jump-noj sumry}~(1) implies that
$s_{N-2}(T_{j_{2n}}) > \tilde{s}_{cr}(T_{j_{2n}}) =
\rho(e_{cr}(T_{j_{2n}}), r)$ (recall that $e_{cr}(T_i)$ is the edge
of $\ch(T_i)$ corresponding to the critical hole in $T_i$).  Thus,
by Proposition~\ref{Prop wsizeto0} and Lemma~\ref{propofrho}~(3),
this implies that $\ch(T_{j_{2n}})$ converges onto $r$ as $n \to
\infty.$ Next, define $w_n = j_{2n} - j_{2n-1}$ for every $n.$  By
Corollary~\ref{crit-iterate}, $f^{w_n}(l)$ lies in one of the $N-2$
smallest holes in $T_{j_{2n}}.$  And by Lemma~\ref{small N-1}(2) we
know that $\lim_{n \to \infty} w_n = \infty.$  These facts together
with Lemma~\ref{trapped seq} imply $\omega(l) \cap r \neq \0.$ Since
$\omega(l) \subset \omega(q)$, this contradicts $r \in \mc{C}_{J}
\setminus C_{\omega(q)}.$
\end{proof}

By the Recurrence Theorem \ref{Recurrence Theorem} there is a unique
$\omega$-limit set for the jumping critical leaves of $T$.  Using
the Recurrence Theorem, we give a specific relationship between the
limit set $\mc{L}_{\omega(T)}$ and the limit set of the jumping
critical leaves.  Note that $\mc{L}_{\omega(T)}$ contains several
chords, while $\omega(q) \subset \ucirc$ for a jumping critical leaf
$q.$  Thus $\mc{L}_{\omega(T)} \neq \omega(q)$, however, we know
that $\omega(q) \subset \mc{L}_{\omega(T)}.$  The following
corollary better describes the relationship between these two limit
sets.

\begin{cor} \label{limcoin}
If $q$ is a jumping critical leaf then for every $p \in
\mc{L}_{\omega(T)}$, $\omega(q) \cap p \neq \0.$
\end{cor}

\begin{proof}
Let $\ch(T_{i_n})$ be a subsequence converging to a leaf $p \in
\mc{L}_{\omega(T)}.$  By Theorem~\ref{Recurrence Theorem}, we need
to show that there exists a jumping critical leaf $q$ with
$\omega(q) \cap p \neq \0.$  To do this, define the subsequence
$\{m_n\}$ of jumps by $m_n = \mathrm{max}\{ i < i_n ~:~ i \in
\jump(T)\}.$  We may assume that $m_{n+1} > m_n$ for every $n$ and
that the critical leaf $q_{m_n}$ (defined in Theorem~\ref{TH
jump-noj sumry}~(1)) is the same critical leaf $q$ for every jump
$m_n.$

Define $w_n = i_n - m_n$. By Corollary \ref{crit-iterate},
$f^{w_n}(q)$ lies in one of the $N-2$ smallest holes in $T_{i_n}$
for every $n$. By Proposition~\ref{trapped seq}, the (non-empty) set
$A$ of limit points of the sequence $f^{w_n}(q)$ is contained in
$p\cap\ucirc.$ Since $q$ is recurrent, $\ol{orb(f(q))} =
\omega(f(q)),$ and so $A\subset \omega(f(q)),$ which completes the
proof.
\end{proof}

\section{Recurrent critical leaves from the limit set of co-existing
wandering polygons} \label{SevPolysSec}

Consider now several co-existing wandering polygons.  Recall that
for a finite set $A \subset \ucirc$, $A_i = f^i(A).$

\begin{dfn}
Let $\Gamma$ be a collection of pairwise disjoint finite subsets.
Then $\Gamma$ is said to \emph{form a wandering collection of
polygons} if the following holds.
\begin{enumerate}

\item
Each finite set $T \in \Gamma$ forms a wandering polygon under $f.$

\item
For every two different finite sets $A, B \in \Gamma$ and any two
integers $n,m \geq 0$, $A_n$ and $B_m$ are unlinked.

\end{enumerate}
\end{dfn}

For a wandering collection of polygons $\Gamma,$ define
$\mc{L}_{\omega(\Gamma)} = \cup_{T \in \Gamma} \mc{L}_{\omega(T)}.$
So, we call any chord (including degenerate chords, i.e., points in
$\ucirc$) in $\mc{L}_{\omega(\Gamma)}$ a \emph{leaf}.  Let $R$ be
the maximum number of recurrent critical leaves with pairwise
disjoint grand orbits. Let $\Omega$ be the collection of distinct
$\omega$-limit sets of the recurrent critical leaves.  Now we can
prove the following theorem.

\begin{thm} \label{wandcol}
Let $f(z) = z^d$ and suppose $\Gamma$ forms a wandering collection
of polygons.  Then $$\card(\Gamma) \leq \sum_{T \in \Gamma}
(\card(T)-2) \leq R - \card(\Omega) \leq d - 1 - \card(\Omega).$$
\end{thm}

\begin{proof} First, let us consider a wandering collection of polygons
$\Gamma'$ with $\card(\Gamma') < \infty.$ Then fix $T \in \Gamma'$
and set $N_T = \card(T)$.  We can construct the sequences $i_n,
j_n(1), \ldots, j_n(N_T-2)$ of jumps as in the proof of
Theorem~\ref{(N-1)critical} applied to this set $T \in \Gamma'$.
Recall that $\tau_n = s_1(T_{i_n + n})$.  As in the proof of
Theorem~\ref{(N-1)critical}, we can find a collection of critical
leaves $\mc{J}_T = \{l_1, \ldots, l_{N_T - 2}\} \subset
\mc{L}_{\omega(\Gamma')}$ and an infinite subset $V \subset
\mathbb{N}$ such that for every $n \in V$, $s_k(T_{j_n(k)}) \geq
\tau_n,$ $q_{j_n(k)} = l_k,$ and $\tau_n > s_k(T_i)$  for all
$i>j_n(k).$  Recall that $l_0 = q_{i_n}$ for every $n.$  As in
Theorem~\ref{(N-1)critical}, $l_0, \ldots , l_{N_T -2}$ are jumping
critical leaves of $T$ with pairwise disjoint orbits.

Let $A \in \Gamma', A \neq T$, and set $N_A = \card(A).$ Since
$\lim_{n \to \infty}\tau_n = 0$ and $V$ is infinite, we can find an
$n \in V$ with  $s_1(A) > \tau_n$. Then by the same methods used to
construct the collection $\mc{J}_T$, we can find a collection of
jumping critical leaves of $A$, $\mc{J}_A = \{l_1', \ldots, l_{N_A
-2}'\} \subset \mc{L}_{\omega(\Gamma')}$, and an infinite subset $V'
\subset V$ such that for all $n \in V'$ there is a subset
$\{j_{n}'(1), j_n'(2), \ldots, j_n'(N_A-2)\} \subset \jump(A)$ such
that $s_k(A_{j_n'(k)}) \geq \tau_n,$ $q_{j_n'(k)} = l_k',$ and
$\tau_n > s_k(A_i)$  for all $i >j_n'(k).$  Since $V' \subset V$,
this implies (as in the proof of Theorem~\ref{(N-1)critical}) the
critical leaves of the collection $\{l_0\} \cup \mc{J}_T \cup
\mc{J}_A\} \subset \mc{L}_{\omega(\Gamma')}$ have pairwise disjoint
orbits. Note that $\card(\mc{J}_A) = \card(A)-2.$

This process of extending our collection of critical leaves can be
continued as long as we have a supply of wandering polygons with
pairwise disjoint orbits. Thus, for every finite wandering
collection of polygons $\Gamma'$ there exists a collection of
critical leaves $\mc{J}_{\Gamma'}= \{l_0\} \cup \cup_{Q \in
\Gamma'}mc{J}_Q \subset \mc{L}_{\omega(\Gamma')}$ with pairwise
disjoint orbits such that $\card(\mc{J}_{\Gamma'}) = \sum_{Q \in
\Gamma'} (\card(Q)-2) + 1,$ and every $l \in \mc{J}_{\Gamma'}$ is a
jumping critical leaf of some $A \in \Gamma'.$  For any finite
wandering collection of polygons $\Gamma',$ let us denote by
$\mc{J}_{\Gamma'}$ the just-constructed collection of critical
leaves; let us emphasize here that all critical leaves from this
collection have pairwise disjoint forward orbits.

Now let $\Gamma$ be the wandering collection of polygons.  We proof
the inequality for $\Gamma.$  First, since the number of critical
leaves is finite, and since $\card(Q)-2 \geq 1$ for any wandering
$N$-gon, by applying the above paragraph to finite subcollections of
$\Gamma,$ we have that $\card(\Gamma) < \infty.$  Now, let us begin.

By the Recurrence Theorem~\ref{Recurrence Theorem}, all of the
jumping critical leaves of a polygon are recurrent and have the same
$\omega$-limit set.  For every $A \in \Gamma$, let us denote by
$W_A$ this limit set.  And for every $X \in \Omega$, define the
subcollection $\Gamma_X = \{ A \in \Gamma : W_A = X\}.$  Set
$\Omega_1 = \{X \in \Omega: \Gamma_X \neq \0\}$ and $\Omega_0 = \{X
\in \Omega: \Gamma_X = \0\}.$ Notice that these two sets are a
partition of $\Omega,$ so $\card(\Omega_0) + \card(\Omega_1) =
\card(\Omega).$

By Theorem~\ref{(N-1)critical} every wandering $N$-gon has at least
$N - 1 > 0$ jumping critical leaves, and so the collection
$\{\Gamma_X\}_{X \in \Omega_1}$ is a partition of $\Gamma$ into
non-empty (by the choice of $\Omega_1$) sets.  For every $X \in
\Omega_1$, we have $\card(\mc{J}_{\Gamma_X})=1+\sum_{Q \in \Gamma_X}
(\card(T)-2)$. Summing up these formulas over all $X\in \Omega_1$,
we see that $\sum_{X\in \Omega_1}\card(\mc{J}_{\Gamma_X}) =
\card(\Omega_1) + \sum_{T \in \Gamma} (\card(T)-2).$

Now, for every $X \in \Omega_0$ there exists at least one recurrent
critical leaf whose limit set is $X$. So we can define a collection
$\mc{R}_0$ of recurrent critical leaves such that for every $X \in
\Omega_0$ there is a unique $l \in \mc{R}_0$ with $\omega(l) = X.$
Note then that $\card(\mc{R}_0) = \card(\Omega_0).$

Finally, let us show that critical leaves in the collection $\cup_{X
\in \Omega_1} mc{J}_{\Gamma_X} \cup \mc{R}_0 $ have pairwise
disjoint forward orbits. Indeed, if the leaves in question come from
the same collection $\mc{J}_{\Gamma_X}$ for some $X$, then it
follows from the properties of such collections. On the other hand,
if they come from different collections $\Gamma_X, \Gamma_Y$ with
$X\ne Y,$ or $\Gamma_X, \mc{R}_0$ for some $X,$ or if they both come
from $\mc{R}_0$, then, in either one of these three cases, they have
different limit sets which again yields that their forward orbits
are disjoint. Thus, $\sum_{X\in \Omega_1}\card(\mc{J}_{\Gamma_X}) +
\card(\mc{R}_0) \le R,$ which implies that $\sum_{T \in \Gamma}
(\card(T)-2) \le R - (\card(\Omega_1)+ \card(\Omega_0)) \le R -
\card(\Omega),$ as desired.  \end{proof}

\bibliographystyle{plain}
\bibliography{e:/lex/references/refshort}

\end{document}